\documentclass[10pt,leqno]{amsart}
\usepackage{graphicx}
\usepackage{indentfirst,csquotes}

\usepackage{amssymb,amsthm,amsmath}
\usepackage{xcolor,paralist,hyperref,fancyhdr,etoolbox}
\newtheorem{th.}{Theorem}[]
\newtheorem{de.}[th.]{Definition}
\newtheorem{ex.}[th.]{Example}
\newtheorem{le.}[th.]{Lemma}
\newtheorem{re.}[th.]{Remark}
\newtheorem{pr.}[th.]{Proposition}
\newtheorem{co.}[th.]{Corollary}
\newtheorem{conj.}[th.]{Conjecture}

\hypersetup{ colorlinks=true, linkcolor=black, filecolor=black, urlcolor=black }

\usepackage{lipsum}

\begin{document}
\title{Quandles from gauge transformations} 
\author{Ryo Hayami}
\date{\today}
\address{Nagano University}
\email{ryo-hayami@nagano.ac.jp}

\let\thefootnote\relax
\footnotemark 

\begin{abstract}
In this paper, we investigate a quandle structure induced by an augmented rack arising from a gauge transformation group. We construct a quandle from a principal bundle and its discrete generalization. When we see a group as a (discrete) principal bundle over a point, this quandle becomes equivalent to the generalized Alexander quandle for its inner automorphism. Moreover, we construct a Lie and Noether quandle structure from a smooth gauge transformation. 
\end{abstract}
\maketitle

\tableofcontents

\bigskip

\section{introduction}

Racks and quandles are algebraic axiomatizations of the Reidemeister moves, which can be regarded as generalizations of the conjugation of groups.

The definition of racks and quandles is summarized as follows.

\begin{de.}

A set $X$ with a binary operation $\triangleleft:X\times X \rightarrow X$ is a rack if it satisfies 
\begin{equation}
(x\triangleleft y)\triangleleft z=(x\triangleleft z)\triangleleft(y\triangleleft z)
\end{equation}
for all $x,y,z\in X$ (this relation is called self-distributivity) and the map $-\triangleleft x:X\rightarrow X$ is bijective. When a rack $(X,\triangleleft)$ satisfies
\begin{equation}
x\triangleleft x=x
\end{equation}
for all $x\in X$, we say that $(X,\triangleleft)$ is a quandle.
When $X$ is a topological space and $\triangleleft$ is continuous, $(X,\triangleleft)$ is called a topological quandle. When $X$ is a smooth manifold and $\triangleleft$ is smooth, $(X,\triangleleft)$ is called a smooth quandle.

For quandles $(X,\triangleleft^{X})$ and $(Y,\triangleleft^{Y})$, a quandle morphism is a map $f:X\rightarrow Y$ such that
\[
f(x_{1}\triangleleft^{X}x_{2})=f(x_{1})\triangleleft^{Y} f(x_{2})
\]
for all $x_{1},x_{2}\in X$. When $f$ is bijective, we say that $f$ is a quandle isomorphism. A quandle morphism between topological quandles which is a homeomorphism and a quandle morphism between smooth quandles which is a diffeomorphism are also called a quandle isomorphism.

\end{de.}

For a group $G$, the conjugation $g_{1}\triangleleft g_{2}=g_{2}^{-1}g_{1}g_{2}$ satisfies the conditions of quandles. This quandle $(G,\triangleleft)$ is called a conjugation quandle. (When we want to focus on the rack structure, we call $(G,\triangleleft)$ a conjugation rack.) When $G$ is a topological (resp. Lie) group, its conjugation quandle becomes a topological (resp. smooth) quandle. 

Although quandles were first introduced in a purely algebraic setting\cite{joyce1982classifying}, topological quandles introduced in \cite{rubinsztein2007topological} and smooth quandles introduced in \cite{Ishikawa} are also important research subjects now. The relation between racks and quandles is summarized in \cite{andruskiewitsch2003racks}. In particular, there is a process for making a quandle $(X,\triangleleft^{\iota})$ from a given rack $(X,\triangleleft)$.

An important class of racks is an augmented rack, whose rack structure arises from the conjugation rack of a group\cite{fenn1992racks}.

\begin{de.}

Let $G$ be a group and act on a set $X$ with a right action. We say that $X$ with a map $\Theta:X\rightarrow G$ is an augmented rack if
\begin{equation}\label{ar}
\Theta(x\cdot g)=g^{-1}\Theta(x)g
\end{equation}
holds for all $g\in G$ and all $x\in X$.

\end{de.}

When $G$ is a topological group and act on a topological space $X$ by homeomorphisms and $\Theta$ is continuous, we say that $(G,X,\Theta)$ is an augmented topological rack. When $G$ is a Lie group and act on a smooth manifold $X$ by diffeomorphisms and $\Theta$ is smooth, we say that $(G,X,\Theta)$ is an augmented smooth rack.

For an augmented rack $(X,G,\Theta)$, $X$ becomes a rack via
\[
x\triangleleft y:=x\cdot\Theta(y).
\]
When the augmented rack is topological (resp.smooth), the induced rack becomes a topological (resp.smooth) rack.

The aim of this paper is to investigate racks and quandles which are given by gauge transformations or their topological or set-theoretical analog.

A gauge transformation is defined as a fiber-preserving automorphism of a principal bundle. Let $\pi:P\rightarrow M$ be a principal bundle over a smooth manifold $M$ whose structure group is a Lie group $G$. For the study of geometry of a gauge theory, we have to comupute the group of the gauge transformations of $P$. The group of gauge transformations is defined as
\[
\mathrm{Gauge}(P):=\{\phi\in \mathrm{Diff}(P)|\forall p\in P,\forall g\in G,  \pi(\phi(p))=\pi(p),\ \phi(p)\cdot g=\phi(p\cdot g)\}.
\]
An important property of the group of gauge transformations is that this group is isomorphic to
\[
C^{\infty}(P,G)^{G}:=\{f\in C^{\infty}(P,G)|f(p\cdot g)=g^{-1}f(p)g\}
\]
with a group structure defined by the pointwise multiplication.

We can see that $(P,G,C^{\infty}(P,G)^{G})$ is an augmented smooth rack. From this, we can define the structure of a smooth rack. With a little more effort, we can associate a smooth quandle structure on $P$ with a gauge transformation on $P$.

The above discussion can be applied even without the smooth setting. Therefore, if we replace the smooth structure with the general topological structure or simply set-theoretical structure, we obtain the topological or set-theoretical analog of quandles from (set-theoretical of topological analog of) gauge transformations. We define and study these quandles, which we call them gauge quandles.

 Treating a group as a principal bundle over a point, we can construct a gauge quandle structure on the group. We show that this quandle is equivalent to a generalized Alexander quandle for an inner automorphism. Moreover, we see that the "fiber quandle" of a gauge quandle is isomorphic to a generalized Alexander quandle for an inner automorphism.

In the smooth setting, by choosing a gauge transformation defining a smooth gauge quandle, we can make an optional quandle-like structure which is called a Lie quandle. A Lie quandle was introduced in \cite{fritz2025self} as a formulation of observables in a physical theory. We show that we can construct a Lie quandle from gauge quandles parametrized by a real number and that this Lie quandle, which we call a parametrized gauge quandle, is a Noether quandle, which is a formulation of Noether's theorem. 

These results would connect the theory of principal bundles and gauge theories with the theory of quandles.  

 \begin{re.}

There is a left self-distributive structure which is expressed as
\begin{equation*}
x\triangleright (y\triangleright z)=(x\triangleright y)\triangleright(x\triangleright z)
\end{equation*}
and in some literatures racks and quandles are defined with this structures. We adopt the right self-distributive structure because a principal bundle is usually defined with a right action. We can easily modify the discussions in this paper to match the left self-distributive structure.

\end{re.}

\section{Quandles from (discrete) gauge transformations}

In this section, we define a rack and a quandle structure from (a discrete or topological generalization of) a gauge transformation.

First, we prepare a discrete version of a principal bundle.

\begin{de.}

Let $P$ and $M$ be sets and let $G$ be a group. We say that $P$ is a discrete principal bundle over $M$ whose structure group is $G$ if there is a surjective map $\pi:P\rightarrow M$ and a right action of $G$ on $P$ 
\[
P\times G\rightarrow P,\ (p,g)\mapsto p\cdot g
\]
which satisfies the conditions.

\begin{itemize}
    \item For each $m\in M$, there is a bijection $\psi_{m}:\pi^{-1}(m)\stackrel{\sim}{\longrightarrow}G$.

    \item For each $p\in \pi^{-1}(m)$, $G$ acts on $p\in \pi^{-1}(m)$ freely and transitively, and they satisfy $\pi(p\cdot g)=\pi(p)=m$ and $\psi_{m}(p\cdot g)=\psi_{m}(p) g$.
\end{itemize}

\end{de.}

When $P$ and $M$ is topological spaces, $G$ is a topological group, $\pi$ is continuous and there is a local trivialization, it becomes an ordinary principal bundle. When $P$ and $M$ is smooth manifolds, $G$ is a Lie group, $\pi$ is smooth and the local trivialization is smooth, it becomes a smooth principal bundle.

Next, we define a discrete and topological generalization of a gauge transformation group.

\begin{de.}

We define a discrete gauge transformation group of a discrete principal bundle as
\[
\mathrm{Gauge}(P)=\{\phi\in \mathrm{Aut}(P)|\pi\circ\phi(p)=\pi(p),\phi(p\cdot g)=\phi(p)\cdot g\}
\]
where $p\in P$ and $g\in G$.

Similarly, we define a topological gauge transformation group of a principal bundle as
\[
\mathrm{Gauge}(P)=\{\phi\in \mathrm{Homeo}(P)|\pi\circ\phi(p)=\pi(p),\phi(p\cdot g)=\phi(p)\cdot g\}
\]
where $p\in P$ and $g\in G$.

\end{de.}

In the smooth setting, the topological gauge transformation group becomes an ordinary gauge transformation group. When there is no risk of confusion, we will simply use the word "gauge transformation group" in the all setting.

We consider the discrete case. Define the group of maps
\[
\mathrm{Map}(P,G)^{G}:=\{f\in \mathrm{Map}(P,G)|f(p\cdot g)=g^{-1}f(p)g\}
\]
whose group structure is defined by the pointwise multiplication. Then, the map
\[
\mathrm{Map}(P,G)^{G}\rightarrow \mathrm{Gauge}(P),\ f\mapsto \phi_{f}
\]
,where $\phi_{f}$ is a gauge transformation defined by $\phi_{f}(p):=p\cdot f(p)$, is an isomorphism of groups. Note that $\phi_{f}^{-1}(p)=p\cdot f(p)^{-1}$.

By the definition of $\mathrm{Map}(P,G)^{G}$, we can see that $(P,G,\mathrm{Map}(P,G)^{G})$ is an augmented rack. Thus a gauge transformation $\phi_{f}$ on $P$ defines a rack structure on $P$ via
\[
p_{1}\triangleleft p_{2}:=p_{1}\cdot f(p_{2}).
\]

Next, we define a quandle structure on $P$. There is a general procedure to make quandles from racks\cite{andruskiewitsch2003racks}. We apply this procedure to augmented racks.

Let $(X,\triangleleft)$ be a rack. Let $\iota:X\rightarrow X$ be a map given by $\iota(x)\triangleleft x=x$ and define a new binary operation $\triangleleft^{\iota}:X\times X\rightarrow X$ by $x\triangleleft^{\iota}y:=\iota(x)\triangleleft y$. It is shown in \cite{andruskiewitsch2003racks} that $(X,\triangleleft^{\iota})$ is a quandle. We say that $(X,\triangleleft^{\iota})$ is the quandle associated to $(X,\triangleleft)$.

Let $(X,G,\Theta)$ be an augmented rack. Then $X$ has a rack product given by
\begin{equation}
x\triangleleft y:=x\cdot \Theta(y).
\end{equation}
for all $x,y\in X$. We consider the quandle associated to an augmented rack.

First, we find the $\iota:X\rightarrow X$ which satifies $\iota(x)\triangleleft x=x$. For an augmented rack, the rack structure is $x\triangleleft y=x\cdot \Theta(y)$. Therefore, if we define $\iota(x):=x\cdot\Theta(x)^{-1}$, we have
\[
\iota(x)\triangleleft x=x\cdot\Theta(x)^{-1}\cdot\Theta(x)=x\cdot e_{G}=x
\]
where $e_{G}$ is a unit of $G$ and we can see that $\iota(x)\triangleleft x=x$. Thus we can conclude that the quandle assocaited to the augmented rack $(X,G,\Theta)$ is $(X,\triangleleft^{\iota})$ where
\[
x\triangleleft^{\iota}y=x\cdot\Theta(x)^{-1}\Theta(y).
\]

We summarize the construction as follows.

\begin{th.}
    
Let $\pi:P\rightarrow M$ be a discrete principal bundle with a structure group $G$ and $f\in Map(P,G)^{G}$ define a gauge transformation $\phi_{f}$ of $P$. Then $P$ has a rack structure via
\begin{align*}
 p_{1}\triangleleft^{f} p_{2}:&=p_{1}\cdot f(p_{1})^{-1}f(p_{2})\\
 &=\phi_{f}^{-1}(p_{1})\cdot f(p_{2}).   
\end{align*}

\end{th.}

We can also check the self-distributivity by direct calculations.

\begin{align*}
(p_{1}\triangleleft^{f} p_{2})\triangleleft^{f} p_{3}&=(p_{1}\cdot f(p_{1})^{-1}f(p_{2}))\triangleleft^{f} p_{3}\\
&=p_{1}\cdot f(p_{1})^{-1}f(p_{2})f(p_{1}\cdot f(p_{1})^{-1}f(p_{2}))^{-1}f(p_{3})\\
&=p_{1}\cdot f(p_{1})^{-1}f(p_{2})f(p_{2})^{-1}f(p_{1})f(p_{1})^{-1}f(p_{1})^{-1}f(p_{2})f(p_{3})\\
&=p_{1}\cdot f(p_{1})^{-1}f(p_{1})^{-1}f(p_{2})f(p_{3})\\
&=p_{1}\cdot f(p_{1})^{-1}f(p_{3})f(p_{3})^{-1}f(p_{1})f(p_{1})^{-1}f(p_{1})^{-1}f(p_{3})\\
&\times f(p_{3})^{-1}f(p_{2})f(p_{2})f(p_{2})^{-1}f(p_{3})\\
&=p_{1}\cdot f(p_{1})^{-1}f(p_{3})\triangleleft^{f} p_{2}\cdot f(p_{2})^{-1}f(p_{3})=(p_{1}\triangleleft^{f} p_{3})\triangleleft^{f}(p_{2}\triangleleft^{f} p_{3}).
\end{align*}

For the topological case, if we replace $\mathrm{Map}(P,G)^{G}$ by the group of continuous maps
\[
C(P,G)^{G}:=\{f\in C(P,G)|f(p\cdot g)=g^{-1}f(p)g\},
\]
the all above discussions can be applied and we obtain a topological quandle on $P$ by the same formula.

Similarly, for the smooth case, if we replace $\mathrm{Map}(P,G)^{G}$ by the group of continuous maps
\[
C^{\infty}(P,G)^{G}:=\{f\in C^{\infty}(P,G)|f(p\cdot g)=g^{-1}f(p)g\},
\]
the discussions can be applied and we obtain a smooth quandle on $P$ by the same formula. We call a quandle of this type a gauge quandle.

\begin{ex.}

We can see every set $X$ as a discrete principal bundle whose structure group is the trivial group $G=\{e\}$. In this case, the possible gauge transformation is only the identity map, and the induced gauge quandle is 
\[
x_{1}\triangleleft^{f} x_{2}=x_{1}\cdot e^{-1}e=x_{1}.
\]
This quandle is called a trivial quandle.

\end{ex.}

The following example says that a gauge quandle over a point is equivalent to a generalized Alexander quandle.

\begin{ex.}

We see a group $G$ as a discrete principal bundle over a point $\pi:G\rightarrow \{pt\}$ whose structure group is $G$ itself. In this case $f\in Map(P,G)^{G}$ is uniquely determined by $f(e)\in G$, where $e$ is the unit of $G$. For a general $g\in G$, the value of $f(g)$ is
\[
f(g)=f(e\cdot g)=g^{-1}f(e)g.
\]
Thus $\mathrm{Im}(f)$ cosists of the conjugacy classes of $f(e)$.

The quandle operation $\triangleleft^{f}:G\times G\rightarrow G$ arising from the augmented rack as a principal bundle over a point is as follows.
\begin{align}
g_{1}\triangleleft^{f} g_{2}&=g_{1}\cdot f(g_{1})^{-1}f(g_{2})\\
&=g_{1}g_{1}^{-1}f(e)^{-1}g_{1}g_{2}^{-1}f(e)g_{2}\\
&=\sigma_{f(e)}(g_{1}g_{2}^{-1})g_{2}.
\end{align}
where $\sigma_{f(e)}$ is the inner automorphism of $G$ for $f(e)\in G$. A quandle of this type is called a generalized Alexander quandle.

\begin{de.}

Let $G$ be a group. For $\sigma\in \mathrm{Aut}(G)$, $G$ has a quandle structure via
\[
g_{1}\triangleleft g_{2}:=\sigma(g_{1}g_{2}^{-1})g_{2}.
\]
This quandle is called a generalized Alexander quandle. 

\end{de.}

When $G$ is a topological (resp. Lie) group and $\sigma\in \mathrm{Homeo}(G)$, (resp. $\sigma\in \mathrm{Diff}(G)$,) the associated generalized Alexander quandle becomes a topological (resp. smooth) quandle.

Comparing the quandle structure from a gauge transformation, we can see that the generalized Alexander quandle for an inner automorphism of $G$ corresponds to a gauge quandle on a principal bundle $\pi:G\rightarrow \{pt\}$.

\end{ex.}

We consider the quandle structure on the "fiber" of a (discrete) principal bundle. Let $(P,\triangleleft^{f})$ be a gauge quandle with the underlying (discrete) principal bundle $\pi:P\rightarrow M$. For $m\in M$, we call $\pi^{-1}(m)$ a fiber of $P$ at $m\in M$. For $p_{1},p_{2}\in\pi^{-1}(m)$, $p_{1}\triangleleft^{f} p_{2}\in\pi^{-1}(m)$ because $\pi(p_{1}\triangleleft^{f} p_{2})=\pi(p_{1})=m$. Thus $\pi^{-1}(m)$ has a well-defined quandle operation as a subquandle of $(P,\triangleleft^{f})$. We call this subquandle $(\pi^{-1}(m),\triangleleft^{f})$ a fiber quandle.

For each $m\in M$, there is a bijection $\psi_{m}:\pi^{-1}(m)\stackrel{\sim}{\longrightarrow}G$. We can endow $G$ with a quandle structure via
\[
g_{1}\triangleleft^{G}g_{2}:=\psi_{m}(\psi_{m}^{-1}(g_{1})\triangleleft^{f}\psi_{m}^{-1} (g_{2})).
\]
It is clear that $\psi_{m}:\pi^{-1}\rightarrow G$ is a quandle isomorphism.

For $f\in\mathrm{Map}(P,G)^{G}$, $f\circ \psi_{m}^{-1}\in\mathrm{Map}(G,G)^{G}$ defines a gauge transformation on $G\rightarrow\{pt\}$. Therefore, $(G,\triangleleft^{G})$ is a gauge quandle and a generalized Alexander quandle for its inner automorphism, which is isomorphic to $(\pi^{-1}(m),\triangleleft^{f})$.

The above discussion leads to the following proposition.

\begin{pr.}

For every gauge quandle $(P,\triangleleft^{f})$ with the structure group $G$, its fiber quandle is always isomorphic to a generalized Alexander quandle of $G$ for its inner automorphism.

\end{pr.}

Next, we study the quandle structure inherited from a reduction of the structure group of a principal bundle.

Let $\pi:P\rightarrow M$ be a discrete principal bundle whose structure group is $G$ and $H$ be a subgroup of $G$. Then we can consider the quotient $P/H$ of $P$ for the group action of $H$ as a subgroup of $G$ on $P$. There is a surjection $\pi_{H}:P/H\rightarrow M$ such that there is a bijection from $G/H$ to $\pi_{H}^{-1}(m)$ for each $m\in M$.

$P$ has the quandle structure derived from a gauge transformation on $P$. We assume that $\mathrm{Im}(f)$ is a subset of the normalizer group of $H$. Then $f(p\cdot H)=Hf(p)=f(p)H$ for all $p\in P$ and, for all $p_{1}\cdot h_{1}\in p_{1}\cdot H$ and $p_{2}\cdot h_{2}\in p_{2}\cdot H$, the quandle operation satisfies

\begin{align}
p_{1}\cdot h_{1}\triangleleft^{f}  p_{2}\cdot h_{2}&=p_{1}\cdot h_{1}f(p_{1}\cdot h_{1})^{-1}f(p_{2}\cdot h_{2})\\
&=p_{1}\cdot f(p_{1})^{-1}h_{1}h_{2}^{-1}f(p_{2})h_{2}\\
&=(p_{1}\cdot f(p_{1})^{-1}f(p_{2}))\cdot \tilde{h} h_{2}\\
&=(p_{1}\triangleleft^{f} p_{2})\cdot\tilde{h} h_{2}\in(p_{1}\triangleleft^{f} p_{2})H
\end{align}
where $\tilde{h}\in H$ such that $h_{1}h_{2}^{-1}f(p_{2})=f(p_{2})\tilde{h}$.

Therefore, $P/H$ has a quandle structure via
\begin{align}
[p_{1}]\triangleleft^{f'} [p_{2}]&=[p_{1}\triangleleft^{f} p_{2}]\\
&=[p_{1}\cdot f(p_{1})^{-1}f(p_{2})]. 
\end{align}
We call the quandle $(P/H,\triangleleft^{f'})$ a reduced quandle of $P$.

\begin{ex.}
Like the previous example, we consider the case for a group $G$ as a discrete principal bundle over a point $\pi:G\rightarrow \{pt\}$ whose structure group is $G$ itself. For a subgroup $H$ of $G$, we assume that every conjugacy class of $f(e)$ is contained in the normalizer of $H$. Then the quandle operation $\triangleleft^{f}:G/H\times G/H\rightarrow G/H$ as the reduced quandle is as follows.
\[
[g_{1}]\triangleleft^{f} [g_{2}]=[\sigma_{f(e)}(g_{1}g_{2}^{-1})g_{2}].
\]
where $\sigma_{f(e)}$ is the inner automorphism of $G$ for $f(e)\in G$.

When we assume that $f(e)h=hf(e)$ for all $h\in H$, we have $\sigma_{f(e)}(h)=h$ for all $h\in H$. We can define a quandle structure on $H\backslash G$ which has the same form as the above case. A quandle of this type is called a homogeneous quandle.

When $H$ is a normal subgroup of $G$, then the two conditions satisfy and the left coset and right coset coincides. Thus the homogeneous quandle of $H\backslash G$ and the reduced quandle of $G/H$ coincides as a quandle. 

\end{ex.}

\section{Smooth gauge transformations and Lie quandles}

In this section, we consider the case of smooth gauge quandles and show that we can make an optional quandle-like structure which is introduced in \cite{fritz2025self} under the name a Lie quandle and Noether quandle. 

Consider a smooth principal bundle $P$ over a manifold $M$ whose structure group is a Lie group $G$. Let $P\times_{G}\mathfrak{g}$ be the adjoint bundle of $P$. For $X\in\Gamma(P\times_{G}\mathfrak{g})$, we have a exponential map
\[
\exp:\Gamma(P\times_{G}\mathfrak{g})\rightarrow\Gamma(P\times_{G}G)\simeq C^{\infty}(P,G)^{G}
\]
arising from the exponential map $\exp:\mathfrak{g}\rightarrow G$. Using this map, for every $t\in\mathbb{R}$, we can define a quandle structure via
\[
p_{1}\triangleleft^{X}_{t} p_{2}:=p_{1}\cdot \exp(-t X(p_{1}))\exp(tX(p_{2}))
\]
where $X\in\Gamma(P\times_{G}\mathfrak{g})$. We call $(P,\triangleleft^{X}_{t})$ a parametrized gauge quandle.

The definition of a Lie quandle and a Noether quandle is following.

\begin{de.}

A Lie quandle is a smooth manifold $P$ with an operation
\[
\triangleleft:P\times\mathbb{R}\times P\rightarrow P
\]
such that the following equation holdsfor all $x,y,z\in P$ and $s,t\in\mathbb{R}$.

\begin{itemize}
    \item self-action 
    \[(x\triangleleft^{X}_{t} y)\triangleleft^{X}_{s}y=x\triangleleft^{X}_{s+t}y.\]
    \item self-distributivity
    \[(x\triangleleft^{X}_{t} y)\triangleleft^{X}_{s}z=(x\triangleleft^{X}_{s}z)\triangleleft^{X}_{t}(y\triangleleft^{X}_{s}z).\]
    \item idempotency
    \[x\triangleleft^{X}_{s}x=x\]
\end{itemize}

When a Lie quandle $(P,\triangleleft)$ satisfies
\[
p_{1}\triangleleft^{X}_{t} p_{2}=p_{1}\ \forall t\in\mathbb{R}\iff p_{2}\triangleleft^{X}_{t} p_{1}=p_{2}\ \forall t\in\mathbb{R},
\]
$(P,\triangleleft)$ is called a Noether quandle.

\end{de.}

A gauge theory originated from the study of theoretical physics, and thus it is an important problem whether gauge quandles have physical meanings. The following theorem would be a hint of this problem.

\begin{th.}

A parametrized gauge quandle $(P,\triangleleft^{X}_{t})$ is a Lie quandle. Moreover, it is a Noether quandle.

\end{th.}

\begin{proof}

Note that
\[
\exp(tX(p_{1}\cdot \exp(sX(p_{2}))))=\exp(-sX(p_{2}))\exp(tX(p_{1}))\exp(sX(p_{2}))
\]
for every $p_{1},p_{2}\in P$.

Let $p_{1},p_{2},p_{3}\in P$ and $s,t\in\mathbb{R}$. Then we can check the axioms as follows.

\bigskip

self-action:

\begin{align*}
(p_{1}\triangleleft^{X}_{t}p_{2})\triangleleft^{X}_{s}p_{2}&=(p_{1}\cdot \exp(-tX(p_{1}))\exp(tX(p_{2})))\triangleleft _{s}p_{2}\\
&=(p_{1}\cdot \exp(-tX(p_{1}))\exp(tX(p_{2})))\exp(-tX(p_{2}))\exp(tX(p_{1}))\exp(-sX(p_{1}))\\
&\ \times \exp(-tX(p_{1}))\exp(tX(p_{2}))\exp(sX(p_{2}))\\
&=p_{1}\cdot\exp(-sX(p_{1}))\exp(-tX(p_{1}))\exp(tX(p_{2}))\exp(sX(p_{2}))\\
&=p_{1}\cdot\exp(-(s+t)X(p_{1}))\exp((s+t)X(p_{2}))=p_{1}\triangleleft^{X}_{s+t}p_{2}.
\end{align*}

\bigskip

self-distributivity:

\begin{align*}
(p_{1}\triangleleft^{X}_{t}p_{2})\triangleleft^{X}_{s}p_{3}&=(p_{1}\cdot \exp(-tX(p_{1}))\exp(tX(p_{2})))\triangleleft _{s}p_{3}\\
&=(p_{1}\cdot \exp(-tX(p_{1}))\exp(tX(p_{2})))\exp(-tX(p_{2}))\exp(tX(p_{1}))\exp(-sX(p_{1}))\\
&\ \times \exp(-tX(p_{1}))\exp(tX(p_{2}))\exp(sX(p_{3}))\\
&=p_{1}\cdot\exp(-sX(p_{1}))\exp(-tX(p_{1}))\exp(tX(p_{2}))\exp(sX(p_{3}))\\
&=p_{1}\cdot\exp(-sX(p_{1}))\exp(sX(p_{3}))\\
&\times\exp(-sX(p_{3}))\exp(sX(p_{1}))\exp(-tX(p_{1}))\exp(-sX(p_{1}))\exp(sX(p_{3}))\\
&\times \exp(-sX(p_{3}))\exp(sX(p_{2}))\exp(tX(p_{2}))\exp(-sX(p_{2}))\exp(sX(p_{3}))\\
&=p_{1}\cdot\exp(-sX(p_{1}))\exp(sX(p_{3}))\triangleleft^{X}_{t} p_{2}\cdot\exp(-sX(p_{2}))\exp(sX(p_{3}))\\
&=(p_{1}\triangleleft^{X}_{s}p_{3})\triangleleft^{X}_{t}(p_{2}\triangleleft^{X}_{s}p_{3}).
\end{align*}

\bigskip

idempotency:

\[
p_{1}\triangleleft^{X}_{s} p_{1}:=p_{1}\cdot \exp(-s X(p_{1}))\exp(sX(p_{1}))=p_{1}.
\]

\bigskip

The condition that $p_{1}\triangleleft^{X}_{t} p_{2}=p_{1}$ for all $t\in\mathbb{R}$ is equivalent to $\exp(-X(p_{1}))\exp(X(p_{2}))=e_{G}$, and $p_{2}\triangleleft^{X}_{t} p_{1}=p_{2}$ for all $t\in\mathbb{R}$. Thus $(P,\triangleleft^{X}_{t})$ is a Noether quandle.

\end{proof}

\section{future directions}

Although our construction can be applied for discrete principal bundle, there are many examples of quandles with topological or smooth structures from topological or smooth principal bundles. 

If we replace $\mathrm{Map}(P,G)^{G}$ with $\mathrm{Map}(P,X)^{G}$ with $(G,X,\Theta)$ be an augmented rack, $P$ has another rack structure called a crossed module of racks\cite{crans2014crossed}. The induced quandle would be worth studying. 

An important structure with the self-distributivity which is not discussed in this paper is a Lie rack, which was introduced in \cite{kinyon2007leibniz} as an integrated object of a Leibniz algebra. The theory around Lie racks may encounter gauge quandles in a generalized gauge theories.

The physical meaning of a parametrized gauge quandle is still ambiguous. Maybe it would be useful for understanding line operators in gauge theories.

$\,$

$\,$

\bibliography{myref}

@article{kinyon2007leibniz,
  title="{Leibniz Algebras, Lie Racks, and Digroups}",
  author={Kinyon, Michael K},
  journal={Journal of Lie Theory},
  volume={17},
  number={1},
  pages={099--114},
  year={2007}
}

@article{fenn1992racks,
  title={Racks and links in codimension two},
  author={Fenn, Roger and Rourke, Colin},
  journal={Journal of Knot theory and its Ramifications},
  volume={1},
  number={04},
  pages={343--406},
  year={1992},
  publisher={World Scientific}
}

@article{andruskiewitsch2003racks,
  title={From racks to pointed Hopf algebras},
  author={Andruskiewitsch, Nicol{\'a}s and Grana, Mat{\i}as},
  journal={Advances in Mathematics},
  volume={178},
  number={2},
  pages={177--243},
  year={2003},
  publisher={Elsevier}
}

@article{rubinsztein2007topological,
  title={Topological quandles and invariants of links},
  author={Rubinsztein, Ryszard L},
  journal={Journal of knot theory and its ramifications},
  volume={16},
  number={06},
  pages={789--808},
  year={2007},
  publisher={World Scientific}
}

@article{fritz2025self,
  title={Self-distributive structures in physics},
  author={Fritz, Tobias},
  journal={International Journal of Theoretical Physics},
  volume={64},
  number={3},
  pages={73},
  year={2025},
  publisher={Springer}
}

@article{joyce1982classifying,
  title={A classifying invariant of knots, the knot quandle},
  author={Joyce, David},
  journal={Journal of Pure and Applied Algebra},
  volume={23},
  number={1},
  pages={37--65},
  year={1982},
  publisher={Elsevier}
}

@unpublished{Ishikawa,
    AUTHOR={Ishikawa, Katsumi},
    TITLE={On the classification of smooth quandles},
    NOTE={preprint}
}

@article{crans2014crossed,
  title={Crossed modules of racks},
  author={Crans, Alissa S and Wagemann, Friedrich},
  journal={Homology, Homotopy and Applications (HHA)},
  volume={16},
  number={2},
  pages={85--106},
  year={2014}
}
\bibliographystyle{unsrt}

\end{document}